\documentclass[11pt]{amsart}
\usepackage[usenames]{color}
\usepackage{amssymb}
\usepackage{graphicx, epsfig}
\usepackage{latexsym, amsfonts, amscd, amsmath}
\usepackage{mathrsfs}
\makeindex 
\input xy
\xyoption{all}

\definecolor{Indigo}{rgb}{0.2,0.1,0.7}
\definecolor{Violet}{rgb}{0.5,0.1,0.7}

\newtheorem{thm}{Theorem}[subsection]
\newtheorem{prop}[thm]{Proposition}
\newtheorem{lem}[thm]{Lemma}
\newtheorem{cor}[thm]{Corollary}

\theoremstyle{definition}

\newtheorem{exa}[thm]{Example}

\theoremstyle{remark}
\newtheorem{rmk}[thm]{Remark}

\numberwithin{equation}{section} \numberwithin{figure}{section}
\numberwithin{table}{section}

\newcommand{\End}{{\operatorname{End}}}

\newcommand{\Hom}{{\operatorname{Hom}}}

\newcommand{\Ima}{{\operatorname{Im}}}

\newcommand{\Jac}{{\operatorname{Jac }}}

\newcommand{\Norm}{{\bf{N }}}

\newcommand{\Spec}{{\operatorname{Spec }}}

\newcommand{\Tr}{{\operatorname{Tr }}}
\newcommand{\val}{{\operatorname{val}}}

\newcommand{\Lie}{{\operatorname{Lie}}}

\newcommand{\gera}{{\frak{a}}}
\newcommand{\gerb}{{\frak{b}}}

\newcommand{\gerj}{{\frak{j}}}

\newcommand{\gerp}{{\frak{p}}}

\newcommand{\gerH}{{\frak{H}}}

\newcommand{\ux}{{\underline{x}}}

\newcommand{\calD}{{\mathcal{D}}}

\newcommand{\calO}{{\mathcal{O}}}

\def\CC{\mathbb{C}}

\def\FF{\mathbb{F}}

\def\HH{\mathbb{H}}

\def\QQ{\mathbb{Q}}
\def\RR{\mathbb{R}}

\def\ZZ{\mathbb{Z}}

\newcommand{\scrA}{{\mathscr{A}}}

\newcommand{\scrC}{{\mathscr{C}}}

\newcommand{\scrL}{{\mathscr{L}}}

\newcommand{\id}{{\noindent}}
\newcommand{\normal}{{\vartriangleleft}}

\newcommand{\arr}{{\; \longrightarrow \;}}

\newcommand{\injects}{{\; \hookrightarrow \;}}

\newcommand{\fbar}{\overline{\mathbb{F}}_p}
\newcommand{\okplus}{{\mathcal{O}_{K^+}}}

\begin{document}
\marginparwidth 50pt
%\setcounter{page}{-1}
%\setcounter{section}{-1}
%\setcounter{subsection}{+1}
%\thispagestyle{empty}
%\renewcommand{\arraystretch}{1.3}
%\begin{tabular}{p{2.4 cm}p{4.5 cm}}
%(Eq1) & $b_0^{p+1} + a_0d_0^p = 0$ \\
%(Eq2) & $b_0a_0^p + a_0c_0^p = 0$ \\
%\end{tabular}
%To write the word RED in red you write:
%\textcolor{Red}{**RED**}\\
%\textcolor{Orange}{**ORANGE**}\\
%\textcolor{Yellow}{**YELLOW**}\\
%\textcolor{Green}{**GREEN**}\\
%\textcolor{Blue}{**Blue**}\\
%\textcolor{Indigo}{**INDIGO**}\\
%\textcolor{Violet}{**Violet**}

\title[Class Invariants]{Class Invariants for Quartic CM Fields}
\author{Eyal Z. Goren \& Kristin E. Lauter}
\address{Department of Mathematics and Statistics,
McGill University, 805 Sherbrooke St. W., Montreal H3A 2K6, QC,
Canada.}\address{Microsoft Research, One Microsoft Way, Redmond,
WA 98052, USA.} \email{goren@math.mcgill.ca;
klauter@microsoft.com} 

\begin{abstract}
One can define class invariants for a quartic primitive CM field
$K$ as special values of certain Siegel (or Hilbert) modular
functions at CM points corresponding to~$K$. Such constructions
were given in \cite{De Shalit Goren} and \cite{Lauter Class
Invariants}. We provide explicit bounds on the primes appearing in
the denominators of these algebraic numbers. This allows us, in
particular, to construct~$S$-units in certain abelian extensions
of~$K$, where~$S$ is effectively determined by~$K$.
\end{abstract}

\maketitle

\section{Introduction}

One of the main problems of algebraic number theory is the
explicit description of ray class fields of a number field~$K$.
Besides the case of the field of rational numbers, the
theory is most advanced in the case where~$K$ is a complex
multiplication (CM) field. Effective constructions are available
using modular functions generalizing the elliptic modular function
$j$; one constructs modular functions as quotients of two modular
forms on a Siegel upper half space and evaluates at CM points
corresponding to~$K$. The values lie in an explicitly determined
extension of the reflex field~$K^\ast$ of~$K$, that depends on the
field over which the Fourier coefficients of the modular function
is defined, on the level of the modular function, and on the
conductor of the order of~$K$ corresponding to the CM point. We
loosely call magnitudes constructed this way \emph{``class
invariants"} of~$K$. The terminology is proposed because when the
Fourier coefficients are rational and the level is~$1$ the values
of the modular function at CM points lie in ray class fields of
$K^\ast$.

An outstanding problem is the effective construction of units in
abelian extensions of number fields, even in the case of complex
multiplication. A solution of this problem is expected to have
significant impact on obtaining additional cases of Stark's
conjecture. The case of cyclotomic units and elliptic units is
well developed, but in higher dimensional cases little  was known.
The essential problem is that divisors of modular functions cannot
be supported at the boundary of the moduli space. The purpose of
this paper is to provide explicit bounds on the primes appearing
in the denominators of class invariants of a primitive quartic CM
field $K$. This yields, in particular, an explicit bound on the
primes dividing the invariants~$u(\gera, \gerb)$ constructed in
\cite{De Shalit Goren}, thus yielding $S$-units lying in a
specific abelian extension of~$K^\ast$ for an explicit finite set
of primes~$S$. It also yields class polynomials for primitive
quartic CM fields whose coefficients are~$S$-integers as
conjectured in \cite{Lauter Class Invariants}.

\section{Elements of small norm in a definite quaternion algebra}
\subsection{Simultaneous embeddings}\label{subsection: simultaneous embeddings}
Let~$B = B_{p, \infty}$ be ``the" quaternion algebra over~$\QQ$
ramified at~$\{p, \infty\}$. Concrete models for~$B$ can be found
in e.g. \cite[p. 98]{Vigneras}. Let~$\Tr$ and~$\Norm$ be the
(reduced) trace and norm on~$B$ and~$x\mapsto \overline{x} =
\Tr(x) - x$ its canonical involution. Let~$R$ be a maximal order
of~$B$. The discriminant of~$R$ is~$p^2$; if we choose
a~$\ZZ$-basis $v_1, \dots, v_4$ for~$R$ then~$\det(\Tr(v_i
\overline{v_j})) = p^2$; cf. \cite[Prop. 1.1]{Pizer}. Further,
using this basis we may identify~$B\otimes \RR$ with~$\RR^4$. The
bilinear form $\langle \alpha, \beta \rangle =
\Tr(\alpha\overline{\beta})$ is represented with respect to this
basis by an integral symmetric~$4\times 4$ matrix~$M$ with even
diagonal entries, which is positive definite and satisfies
$\det(M) = p^2$. It defines an inner product on~$\RR^4$. We let
$\| r\| = \sqrt{\langle r, r\rangle} = \sqrt{2\Norm(r)}$. Note
that the co-volume of~$R$ (the absolute value of the volume of a
fundamental parallelepiped) is~$p$.

Let~$K_i = \QQ(\sqrt{D_i} )\subset B$ be quadratic imaginary
fields,~$D_i$ a square free integer, and let~$\calO_i$ be orders
of~$K_i$ of conductor~$m_i$ contained in~$R$. We assume that one
of the following equivalent conditions holds: (i)~$K_1 \neq K_2$;
(ii)~$K_1$ does not commute with~$K_2$; (iii)~$K_1 \cap K_2 =
\QQ$. Let~$k_i$ be an element of~$\calO_i$ such that~$\{ 1, k_i\}
$ is a basis of~$\calO_i$ and~$k_i$ has minimal possible norm.

Consider the sublattice~$L$ of~$R$ spanned by~$1, k_1, k_2,
k_1k_2$. It is a full-rank sublattice. Hence, $\text{co-vol}(L)
\geq \text{co-vol}(R)$, while on the other hand
\begin{equation}\label{eqn: simlultaneous embeddings}
\text{co-vol}(L)\leq  \| 1\|\cdot \|k_1\|\cdot \| k_2 \|\cdot
\|k_1k_2\| =  4\cdot \Norm(k_1)\cdot \Norm(k_2). \end{equation} We
therefore get
\begin{lem}\label{lemma:simultaneous embeddings special case} Let~$K_i$, $i = 1, 2$, be quadratic
imaginary fields of discriminant~$d_{K_i}$ contained in~$B$ and
let~$\calO_i$ be the order of conductor~$m_i$ of~$K_i$. Assume
that both~$\calO_1, \calO_2$ are contained in~$R$, a maximal order
of~$B$, and that~$K_1 \neq K_2$. Then
\begin{equation} p \leq \frac{(m_1^2d_{K_1} - 1)(m_2^2d_{K_2} - 1)}{4}.
\end{equation}
\end{lem}
\begin{proof}
One verifies that~$k_i = \pm m_i\sqrt{D_i}$ if~$D_i \equiv 2,
3\pmod{4}$, $(\pm 1 \pm m_i\sqrt{D_i})/2$ if~$D_i \equiv 1
\pmod{4}$ and~$m_i$ is odd, and $\pm m_i\sqrt{D_i}/2$ if~$D_i
\equiv 1 \pmod{4}$ and~$m_i$ is even. The norms are, respectively,
$m_i^2\vert D_i \vert, (1 + m_i^2\vert D_i \vert)/4$ and
$m_i^2\vert D_i \vert/4$.
\end{proof}

\subsection{Corollaries} We draw some corollaries from Lemma~\ref{lemma:simultaneous embeddings special
case}. While some of the corollaries are weaker than what can be
drawn from the work of Gross-Zagier \cite{Gross Zagier} and Dorman
\cite{Dorman}, the techniques are much easier and generalize to
higher dimensional situations. Let~$R$ be a maximal order of~$B$,
the quaternion algebra over~$\QQ$ ramified at~$p$ and~$\infty$.
\begin{cor} If~$x, y \in R$
and~$\Norm(x), \Norm(y) < \sqrt{p}/8$ ($< \sqrt{p}/2$ if~$\Tr(x) =
\Tr(y) = 0$) then~$xy = yx$.
\end{cor}
\begin{proof} One reduces to the case of trace zero elements by
replacing~$x$ and $y$ by~$2x - \Tr(x)$ and $2y - \Tr(y)$. If~$x$ and $y$ have
trace zero we get an embedding of the imaginary quadratic orders
$\ZZ[\sqrt{-\Norm(x)}]$,
$\ZZ[\sqrt{-\Norm(y)}]$ into~$R$. The
Corollary follows from the same argument as above, see
Equation~(\ref{eqn: simlultaneous embeddings}), taking~$k_1 = x,
k_2 = y$.
\end{proof}
\begin{cor}\label{cor: elements of small order lie in a quad subring}
There is an order $\calO_1$ of a quadratic imaginary
field~$K_1$,~$\calO_1 \subset R$ such that all elements of~$R$ of
norm less than~$\sqrt{p}/8$ belong to $\calO_1$. In particular,
for every constant~$A < \sqrt{p}/8$ the number of elements of~$R$
of norm~$A$ is at most the number of elements in~$O_{K_1}$ of
norm~$A$, which is at most~$4\sqrt{A}+2$.
\end{cor}
\begin{cor} Let~$j_i$,~$i = 1, 2$, be two
singular~$j$-invariants,~$j_i$ corresponding to an elliptic curve
$E_i$ with complex multiplication by an order~$\calO_i$ of
conductor~$m_i$ in a quadratic imaginary field~$K_i \subset \CC$.
Suppose that~$K_1 \neq K_2$. Let~$\gerp$ be a prime of
$\overline{\QQ}$ dividing~$p$. If~$(j_1 - j_2) \in \gerp$ then
$p\leq \frac{(m_1^2d_{K_1} - 1)(m_2^2d_{K_2} - 1)}{4}$.
\end{cor}
\begin{proof} If~$(j_1 - j_2) \in \gerp$ then~$E_1 \cong E_2
\pmod{\gerp}$. Let~$E = E_1\pmod{\gerp}$. Since~$K_1 \neq K_2$,
$E$ is supersingular and after fixing an isomorphism
$\End(E)\otimes \QQ \cong B$,~$\End(E)$ is a maximal order of~$B$
containing~$\calO_1, \calO_2$.
\end{proof}

\subsection{A remark on successive minima}
Let~$E, E'$ be two supersingular elliptic curves over
$\overline{\FF}_p$. Then~$\Hom(E, E')$ is a free abelian group of
rank~$4$ equipped with a quadratic form~$f\mapsto \deg(f)$. The
associated bilinear form is~$\langle f, g\rangle \mapsto \deg(f+g)
- \deg(f) - \deg(g) =  f\circ g^\vee + g\circ f^\vee$. Let~$M$ be
a matrix representing~$\langle \cdot, \cdot \rangle$. It is known
to be an integral symmetric positive definite matrix whose
diagonal entries are even and whose determinant is~$p^2$; cf.
\cite[Prop. 1.1]{Pizer}. We are interested in studying the
successive minima~$\mu_1, \dots, \mu_4$ of the homogenous function
of weight~$1$ (``gauge function"), $f(x) = (\frac{1}{2}^tx M
x)^{1/2}$. Geometry of numbers, see \cite[III \S\S3-4, X
\S3]{Siegel}, gives~$2^4(4! V)^{-1}\leq \mu_1 \cdot \mu_2 \cdot
\mu_3 \cdot \mu_4 \leq 2^4V^{-1}$, where $V$ is the volume of the
unit ball with respect to $f$. Since $V = 2\pi^2/p$ we find that
\begin{equation}\label{eqn: successive minima}\frac{1}{3\pi^2}\cdot p\leq \mu_1 \cdot \mu_2 \cdot \mu_3 \cdot
\mu_4 \leq \frac{8}{\pi^2}\cdot p.\end{equation}

\begin{prop} Assume that~$E = E'$. Let~$x$ be an element for which~$\mu_2$ is
obtained,~$K_1 = \QQ(x)$. Then~$\ZZ[x]$ is an order of~$K_1$,
optimally embedded in~$\End(E)$. We have $\mu_2^2 \leq
\frac{4}{\pi^{4/3}}\cdot p^{2/3}$ and
\[\frac{1}{2}\cdot p^{1/2} \leq \mu_2\mu_3 \leq
\frac{4}{\pi^{4/3}}\cdot p^{2/3}, \quad
\max\left\{\frac{1}{3^{1/3}\cdot \pi^{2/3}}\cdot p^{1/3},
\frac{1}{\sqrt{2}}\cdot p^{1/4}\right\}\leq \mu_4\leq
\frac{2\sqrt{2}}{\pi}\cdot p^{1/2}.\]
\end{prop}
\begin{proof} If~$E = E'$, then~$\mu_1 = 1$. The embedding is optimal because the element of
minimal norm independent of~$1$ in an order determines the order,
cf. proof of Lemma~\ref{lemma:simultaneous embeddings special
case}. Let~$y$ be an element for which the third successive
minimum is obtained and~$K_2 = \QQ(y)$. By definition~$y\not\in
K_1$ and we are in the situation of \S\ref{subsection:
simultaneous embeddings}. We get using Equation~(\ref{eqn:
simlultaneous embeddings}), $\sqrt{p} \leq
2\sqrt{\Norm(x)}\sqrt{\Norm(y)} = 2\mu_2\mu_3$; since $xy$ is
independent of $\{1, x, y\}$ we deduce that $\mu_4\leq\sqrt{N(xy)}
= \mu_2\mu_3$. We also have $\mu_4^2\geq \mu_2\mu_3\geq 1$.
Analysis of the inequalities gives the result.
\end{proof}
By Equation~(\ref{eqn: successive minima}) $\mu_1^2 \leq
\frac{2\sqrt{2}}{\pi}\sqrt{p}$; applying that to any pair of
elliptic curves, we get:
\begin{prop} Fix a supersingular elliptic curve~$E$ over
$\overline{\FF}_p$. Let~$N$ be the minimal integer such that for
any supersingular~$E'$ over~$\overline{\FF}_p$ there is an isogeny
of degree less or equal to~$N$ between~$E$ and~$E'$. Then
\[ N \leq \frac{2\sqrt{2}}{\pi}\sqrt{p} \sim 0.9004\sqrt{p}.\]
\end{prop}
\begin{rmk} By estimating the number of subgroups of~$E$ of order
$\leq N$, one finds that~$N \geq 0.3 \sqrt{p}$. Numerical evidence
shows that~$N \geq 0.7 \sqrt{p}$. See~\S\ref{data:isogenies}. Such
results can be interpreted in terms of coefficients of theta
series.
\end{rmk}
\section{An embedding problem}\label{sec: embedding problem}
Let~$K$ be a primitive quartic CM field and~$K^+$ its totally real
subfield. Write~$K^+ = \QQ(\sqrt{d})$, for~$d>0$ a square free
integer. Write~$K = K^+(\sqrt{r})$ with~$r \in \ZZ[\sqrt{d}]$ a
totally negative element. Every quartic CM field can be written
this way and there is much known on the index of~$\ZZ[\sqrt{d},
\sqrt{r}]$ in~$\calO_K$. The following are equivalent for CM
fields of degree~$4$: (i)~$K$ is primitive, i.e., does not contain
a quadratic imaginary field; (ii) $K$ is either non-Galois, or a
cyclic Galois extension; (iii) $\Norm_{\QQ(\sqrt{d})/\QQ}(r)$ is
not a square in $\QQ$.

Let~$E_1, E_2$ be supersingular elliptic curves over
$\overline{\FF}_p$. Let~$\gera = \Hom(E_2, E_1)$,~$\gera^\vee:=
\Hom(E_1, E_2)$, $R_i =\End(E_i)$. Then~$\End(E_1\times E_2) =
\left(\begin{smallmatrix} R_1 & \gera \\ \gera^\vee &
R_2\end{smallmatrix}\right)$. The product polarization induced by
the divisor~$E_1\times \{ 0 \} + \{ 0 \} \times E_2$ on~$E_1
\times E_2$ induces a Rosati involution denoted by $\vee$. This
involution is given by~$\left(\begin{smallmatrix} a & b \\ c & d
\end{smallmatrix}\right) \mapsto \left(\begin{smallmatrix} a & b \\ c
& d \end{smallmatrix}\right)^\vee = \left(\begin{smallmatrix}
a^\vee & c^\vee \\ b^\vee & d^\vee
\end{smallmatrix}\right)$, where~$a^\vee, b^\vee$ etc. denotes the dual
isogeny. Note that~$a^\vee = \overline{a}, d^\vee = \overline{d}$.
The Rosati involution is a positive involution.

\

\id {\bf The embedding problem:} \emph{To find a ring embedding
$\iota\colon \calO_K \injects \End(E_1\times E_2)$ such that the
Rosati involution coming from the product polarization induces
complex conjugation on~$\calO_K$.}

\

\id As we shall see below, the problem is intimately related with
bounding primes in the denominators of class invariants.

\begin{thm}\label{thm: bound for solving the embedding problem}
If the embedding problem has a positive solution then~$p\leq
16\cdot  d^2 (\Tr(r))^2$.
\end{thm}
\begin{proof} Assume such an embedding~$\iota$ exists. Then~$\iota(\calO_{K^+})$ is fixed by the Rosati
involution, thus~$\sqrt{d} \mapsto M = \left(\begin{smallmatrix} a
& b \\ c & e \end{smallmatrix}\right)$, for some~$a, e \in \ZZ, \;
b\in \gera, \; b^\vee = c$. Moreover,~$M^2 =
\left(\begin{smallmatrix} d & 0 \\ 0 & d
\end{smallmatrix}\right)$. This gives the following conditions on the entries of~$M$.

\

{\renewcommand{\arraystretch}{1.3}
   \begin{tabular}{p{1cm}p{4cm}p{4cm}}
 &$a^2 + b{b^\vee} = d$ &$ b(a+e) = 0$ \\
 &${b^\vee}(a+e) = 0$ &${b^\vee} b + e^2 = d$.
   \end{tabular}}

\

\id If~$a \neq -e$ then~$b = 0$ and hence~$d$ is a square - a
contradiction. Thus,~$a = -e$, and we can write the embedding as
\begin{equation} \label{982734}\sqrt{d} \mapsto
  M = \begin{pmatrix}
    a & b \\
    b^\vee & -a
  \end{pmatrix}, \quad a\in \ZZ, \; b\in \gera, \; a^2 + bb^\vee =
  d.\end{equation}
  Let us write
  \[r = \alpha +
  \beta\sqrt{d}, \quad \alpha < 0 , \vert \alpha \vert > \vert
  \beta\sqrt{d} \vert.\]
The condition of the Rosati involution
  inducing complex conjugation is equivalent to~${\iota(\sqrt{r})^\vee}= - \iota(\sqrt{r}).$
  So,
  if~$\iota(\sqrt{r})=
  \left(\begin{smallmatrix}
    x & y \\
    z & w
  \end{smallmatrix}\right)$ then~$\left(\begin{smallmatrix}
    x^\vee & z^\vee \\
    y^\vee & w^\vee
  \end{smallmatrix}\right) =
  -\left(\begin{smallmatrix}
    x & y \\
    z & w
  \end{smallmatrix}\right).$
  This translates into the conditions~$x=-x^\vee, \quad w=-w^\vee, \quad y=-z^\vee,$ implying in
particular that~$x$ and~$w$ have trace zero, which we write as~$x
\in R_1^0$ and~$w \in R_2^0$.  It follows that
\begin{equation} \iota(\sqrt{r}) =  \begin{pmatrix}
    x & y \\
 -y^\vee & w
  \end{pmatrix}, \qquad x \in R_1^0, \; w \in R_2^0, \; y \in \gera.
  \end{equation}
A further condition is obtained from~$\iota(\sqrt{r})^2 = r$,
i.e.,~$\left(\begin{smallmatrix} x & y \\ -y^\vee & w
\end{smallmatrix}\right)^2  = \left(\begin{smallmatrix} \alpha + \beta a & \beta b \\
\beta b^\vee & \alpha - \beta a \end{smallmatrix}\right)$, that is,
\begin{equation}
\begin{pmatrix}
   x^2-yy^\vee & xy+yw \\
   -y^\vee x-wy^\vee & w^2-y^\vee y
\end{pmatrix} =
\begin{pmatrix}
   \alpha + \beta a & \beta b \\
   \beta b^\vee & \alpha - \beta a
\end{pmatrix}.
\end{equation}
Since~$yy^\vee = y^\vee y \in \ZZ$, this leads to the following
necessary conditions

\

  {\renewcommand{\arraystretch}{1.3}
   \begin{tabular}{p{1cm}p{4cm}p{4cm}}
   ($\star$) &$x^2-yy^\vee=\alpha+\beta a$ &$ xy+yw=\beta b$ \\
   &$w^2-yy^\vee = \alpha-\beta a$ &$a^2+bb^\vee = d,$
   \end{tabular}}

\

\id where~$x\in R^0_1, \; w\in R^0_2, \; b,y \in \gera, \;
   \alpha, \beta, a \in \ZZ$.

Note that~$y = 0$ implies that either~$b = 0$ or~$\beta = 0$. The
case~$b = 0$ gives that~$d$ is a square, hence is not possible;
the case~$\beta = 0$ is possible, but leads to~$K$ a bi-quadratic
field, contrary to our assumption.

We use the notation~$\Norm(s) = ss^\vee, \Norm (y) = yy^\vee$,
etc.. Note that for~$s\in R_i$ this definition of the norm is the
usual one and, in any case, under the interpretation of elements
as endomorphisms~$\Norm(s) = \deg(s)$ and so $\Norm(st)
=\Norm(s)\cdot \Norm(t)$ when it makes sense. It follows
from~$(\star)$ that

\

{\renewcommand{\arraystretch}{1.3}
\begin{tabular}{p{1cm}p{7cm}}
$(\star\star)$ &$\Norm(x) + \Norm(y) =-(\alpha+\beta a)$  \\
&$\Norm(w) + \Norm(y) =-(\alpha-\beta a)$.
\end{tabular}}

\

\id Let~$\varphi\colon E_1\rightarrow E_2$ be a non-zero isogeny
of degree~$\delta$. For~$f\in \End(E_1\times E_2)$ the composition
of rational isogenies
\[ E_1 \times E_1 \overset{(1, \varphi)}{\arr} E_1 \times E_2 \overset{f}{\arr} E_1 \times E_2
 \overset{(1, \delta^{-1}\varphi^\vee)}{\arr} E_1 \times E_1,\]
gives a ring homomorphism~$\End^0(E_1 \times E_2) \arr \End^0(E_1
\times E_1)$ that can be written in matrix form as
\[ f = \begin{pmatrix}
  f_{11} & f_{12} \\
  f_{21} & f_{22}
\end{pmatrix} \mapsto \begin{pmatrix}
  f_{11} & f_{12}\varphi \\
  \delta^{-1}\varphi^\vee f_{21} & \delta^{-1}\varphi^\vee f_{22} \varphi
\end{pmatrix}.\]
Let~$\psi$ be the composition~$K \arr \End^0(E_1 \times E_2) \arr
\End^0(E_1 \times E_1).$ Then~$\psi$ is an embedding of rings with
the property
\begin{equation}
  \begin{pmatrix}
    1 & 0 \\
    0 & \delta
  \end{pmatrix}
 \psi(\calO_K) \subset M_2(R_1).
 \end{equation}
Choose~$\varphi = y^\vee$. Taking~$f =  \left(\begin{smallmatrix}
  a & b \\
  b^\vee & -a
\end{smallmatrix}\right)$ (corresponding to~$\sqrt{d}$), or~$\left(\begin{smallmatrix}
  x & y \\
  -y^\vee & w
\end{smallmatrix}\right)$ (corresponding to~$\sqrt{r}$), the embedding~$\psi$ is determined by
\begin{equation} \psi(\sqrt{d}) =
  \begin{pmatrix}
    a & by^\vee \\
    yb^\vee/\delta & -a
  \end{pmatrix}, \qquad
  \psi(\sqrt{r}) =
  \begin{pmatrix}
    x & \delta \\
    -1 & ywy^\vee/\delta
  \end{pmatrix}.\end{equation}
We conclude that
\[ S = \{ by^\vee,y b^\vee, x, y w y^\vee \}
\subset R_1.\] Let

\

{\renewcommand{\arraystretch}{1.3}
\begin{tabular}{lll}$\delta_1$ &~$= \min\{
-(\alpha-\beta a), -(\alpha+\beta a)\}$&$= \vert \alpha \vert -
\vert \beta \vert \cdot \vert a \vert,$\\
$\delta_2$ &  =~$\max\{ -(\alpha-\beta a), -(\alpha+\beta a)
\}$&$= \vert \alpha \vert + \vert \beta \vert \cdot \vert a
\vert.$
\end{tabular}}

\

\id It follows from~(\ref{982734}) and~$(\star\star)$
\[\Norm( by^\vee ) = \Norm( y b^\vee) \leq d\delta_1, \quad \Norm(
x) \leq \delta_2.\] Assume that~$p >  16\cdot d^2(\Tr(r))^2 \ge \max\{ 64 d^2\delta_1^2, 64
\delta_2^2 \}.$ Then $\Norm( by^\vee )$, $\Norm( y
b^\vee)$ and ~$\Norm( x)$ are all smaller than~$\sqrt{p}/8$. By
Corollary~\ref{cor: elements of small order lie in a quad
subring}, the elements~$x, y^\vee b, y b^\vee$ belong to some
imaginary quadratic field $K_1$. The equation~$xy + yw = \beta b$
appearing in~$(\star)$ gives the relation~$xyy^\vee + ywy^\vee =
\beta b y^\vee$, which shows that~$ywy^\vee\in K_1$. We conclude
that~$\psi$ is an embedding~$K \rightarrow M_2(K_1)$. This implies
that~$K_1\injects K$ (else consider the commutative subalgebra
generated by~$K$ and~$K_1$ in~$M_2(K_1)$), contrary to our
assumption. It follows that if there is a solution to the
embedding problem~$\iota$ then~$p \leq 16\cdot
d^2(\Tr(r))^2.$
\end{proof}

\section{Bad reduction of CM curves} In this section we discuss
the connection between solutions to the embedding problem and bad
reduction of curves of genus two whose Jacobian has complex
multiplication. We shall assume CM by the full ring of integers,
but the arguments can easily be adapted to CM by an order, at
least if avoiding primes dividing the conductor of the order.

\subsection{Bad reduction solves the embedding problem}
Fix a quartic primitive CM field~$K$. Write~$K =
\QQ(\sqrt{d})(\sqrt{r})$,~$r\in \calO_{K^+}$,~$d$ a positive
integer.  Let~$\scrC$ be a smooth projective genus~$2$ curve over
a number field~$L$. We say that~$\scrC$ has CM (by~$\calO_K$)
if~$\Jac(\scrC)$ has CM by~$\calO_K$. By passing to a finite
extension of~$L$ we may assume that~$\scrC$ has a stable model
over~$\calO_L$ and that all the endomorphisms of~$\Jac(\scrC)$ are
defined over~$L$. Since $K$ is primitive,~$\Jac(\scrC)$ is a
simple abelian variety and so $\End^0(\Jac(\scrC)) = K$. In
particular, the natural polarization of~$\Jac(\scrC)$, associated
to the theta divisor~$\scrC \subset \Jac(\scrC)$, preserves the
field~$K$ and acts on it by complex conjugation.

It is well known that~$\Jac(\scrC)$ has everywhere good reduction.
It follows that for every prime ideal~$\gerp\normal \calO_L$
either~$\scrC$ has good reduction modulo~$\gerp$ or is
geometrically isomorphic to two elliptic curves~$E_1, E_2$
crossing transversely at their origins. In the latter case we have
an isomorphism of principally polarized abelian varieties
over~$k(\gerp) = \calO_L/\gerp$ , $(\Jac(\scrC), \scrC) \cong (E_1
\times E_2, E_1\times \{ 0 \} + \{ 0 \} \times E_2)$. Since
$K\injects \End(E_1\times E_2)\otimes \QQ$ we see that $E_1$ must
be isogenous to $E_2$. Moreover, $E_i$ cannot be ordinary; that
implies that $K\injects M_2(K_1)$ for some quadratic imaginary
field $K_1$ and one concludes that $K_1 \injects K$, contradicting
the primitivity of~$K$. We conclude

\begin{lem} Let~$\scrC/L$ be a non-singular projective curve of
genus~$2$ with CM by~$\calO_K$. Assume that~$\scrC$ has a stable
model over~$\calO_L$. If~$\scrC$ has bad reduction modulo a prime
$\gerp\vert p$ of~$\calO_L$ then the embedding problem has a
positive solution for the prime~$p$.
\end{lem}
The following theorem now follows immediately using
Theorem~\ref{thm: bound for solving the embedding problem}.
\begin{thm}  Let~$\scrC$ be a non-singular projective curve of genus~$2$ with CM
by~$\calO_K$ and with a stable model over the ring of
integers~$\calO_L$ of some number field~$L$. Let~$\gerp \vert p$
be a prime ideal of~$\calO_L$. Assume that~$p$ is greater or equal
to~$16\cdot  d^2 (\Tr(r))^2$ then~$\scrC$ has good reduction
modulo~$\gerp$.
\end{thm}

\subsection{A solution to the embedding problem implies bad reduction}
\begin{thm}\label{thm: embedding implies bad reduction} Assume that the embedding problem of \S\ref{sec: embedding
problem} has a solution with respect to a primitive quartic CM
field~$K$. Then there is a smooth projective curve~$\scrC$ of
genus~$2$ over a number field~$L$ with CM by~$\calO_K$, whose
endomorphisms and stable model are defined over~$\calO_L$, and a
prime~$\gerp$ of~$\calO_L$ such that~$\scrC$ has bad reduction
modulo~$\gerp$.
\end{thm}
\id Our strategy for proving the theorem is the following. We
consider a certain infinitesimal deformation functor~${\bf N}$ for
abelian surfaces with CM by~$\calO_K$. We show that~${\bf N}$ is
pro-representable by a $W(\overline{\FF}_p)$-algebra~$R^{\rm u}$,
and that a solution to the embedding problem can be viewed as an
$\overline{\FF}_p$-point~$x$ of~$\Spec(R^{\rm u})$. We prove
that~$R^{\rm u}$ is isomorphic to the completed local ring of a
point on a suitable Grassmann variety and deduce that~$R^{\rm u}
\otimes \QQ \neq 0$. We conclude that~$x$ can be lifted to
characteristic zero and finish using classical results in the
theory of complex multiplication. Before beginning the proof
proper, we need some preliminaries about Grassmann varieties.
\subsection{Grassmann schemes}\label{sec: Grassmannians}
\id The following applies to any number field~$K$ with an
involution~$\ast$; we denote the fixed field of~$\ast$ by~$K^+$.
Put $[K:\QQ] = 2g$.

\subsubsection{} Consider the module~$M_0:=\calO_K
\otimes_\ZZ W$,~$W = W(\fbar)$, equipped with an alternating
perfect~$W$-linear pairing~$\langle \cdot, \cdot \rangle$ with
values in~$W$, such that for~$s\in \calO_K$ we have ~$\langle sr,
r' \rangle = \langle r, s^\ast r' \rangle$. Note that this also
holds for~$s\in \calO_K\otimes_\ZZ W$ if~$\ast$ denotes the
natural extension of the involution to this ring.

This defines a Grassmann problem: classify for~$W$-algebras~$W'$
the isotropic, locally free, locally direct summands
$W'$-submodules of~$M_0\otimes_W W'$ of rank~$g$ that
are~$\calO_K$-invariant. This is representable by a projective
scheme~${\bf G'} \rightarrow \Spec(W)$ (a closed subscheme of the
usual (projective) Grassmann scheme). We claim that~${\bf G'}$ is
topologically flat: namely, that every~$\fbar$-point of it lifts
to characteristic zero. That means that for every submodule~$N_1$
of~$\calO_K \otimes_\ZZ\fbar$, satisfying the conditions above,
there is a flat~$W$-algebra~$W'$ and such submodule~$N_0$
of~$\calO_K \otimes_W W'$ that lifts~$N_1$.

\subsubsection{}\label{subsubsection:3765} First note that for~$k\supset W$ an algebraically closed field
of characteristic zero, the~$k$-points of~${\bf G'}$ are in
bijection with ``CM types". Indeed, we are to classify the
isotropic, rank~$g$, sub~$k$-vector spaces of~$\calO_K \otimes_\ZZ
k = \oplus_{\{\varphi\colon K \rightarrow k\}} k(\varphi)$,
where~$k(\varphi)$ is~$k$ on which~$\calO_K$ acts via~$\varphi$.
It is easy to see that the pairing decomposes as a direct sum of
orthogonal pairings on the~$g$ subspaces~$k(\varphi)\oplus
k(\varphi\circ \ast)$ (use that for~$r\in k(\varphi), r' \in
k(\varphi')$ we have~$\varphi(s) \langle r, r'\rangle = \langle
sr, r'\rangle= \langle r, s^\ast r'\rangle =
(\varphi'\circ\ast)(s) \langle r, r'\rangle$). On
$k(\varphi)\oplus k(\varphi\circ \ast)$ the pairing is
non-degenerate so every maximal isotropic subspace is a line and
vice-versa. The condition of being an $\calO_K$-submodule leaves
us with precisely two submodules of~$k(\varphi)\oplus
k(\varphi\circ \ast)$, viz.~$k(\varphi)$, $k(\varphi\circ \ast)$.
Thus, the choice of an isotropic, $\calO_K$-invariant $k$-subspace
of dimension~$g$ of $\calO_K \otimes_\ZZ k$ corresponds to
choosing an element from each of the~$g$ pairs~$\{ \varphi,
\varphi\circ \ast\}$.

\subsubsection{}  We now prove topological flatness for ${\bf G'}$.
We first make a series of reductions. Let~$p = \prod
\gerp_i^{e_i}$ in~$\okplus$. We have $\okplus\otimes_\ZZ W =
\oplus_{\gerp\vert p} \calO_{K^+_\gerp} \otimes_{\ZZ_p} W$ (with
corresponding idempotents~$\{e_\gerp\}$) and~$\okplus\otimes_\ZZ
\fbar = \oplus_{\gerp\vert p} \okplus/\gerp^{e_\gerp}
\otimes_{\FF_p} \fbar$. The modules~$M_0 = \calO_K \otimes_\ZZ
W$,~$M_1 := \calO_K \otimes_\ZZ \fbar$, which are, respectively,
free~$\okplus\otimes_\ZZ W$ and~$\okplus\otimes_\ZZ \fbar$ modules
of rank~$2$, decompose accordingly as~$\oplus_{\gerp\vert p}
M_0(\gerp), \oplus_{\gerp\vert p} M_1(\gerp)$. We claim that the
submodules~$\{ M_0(\gerp) : \gerp\vert p\}$ (resp.~$\{ M_1(\gerp)
: \gerp\vert p\}$) are orthogonal. Indeed, this follows from the
fact that for the idempotents~$\{ e_\gerp\}$ we have~$\langle
e_\gerp r, e_{\gerp'}r' \rangle = \langle e_\gerp^2 r,
e_{\gerp'}^2r' \rangle =\langle e_\gerp r, e_\gerp e_{\gerp'}^2 r'
\rangle = \delta_{\gerp, \gerp'}\langle e_\gerp r, e_{\gerp'} r'
\rangle~$. We may thus assume without loss of generality that~$p =
\gerp^{e}$ with residue degree~$f$ in~$\okplus$ (note that the
global nature of the rings~$\calO_K, \okplus$ plays no role).
Let~$W^+ = W(\FF_{p^f})$ considered as the maximal unramified
sub-extension of~$\calO_{K^+_\gerp}$. A further reduction is
possible: Since~$\calO_{K^+_\gerp} \otimes_{\ZZ_p} W =
\oplus_{\{W^+ \rightarrow W\}}\calO_{K^+_\gerp} \otimes_{W^+} W$,
the same arguments as above (using idempotents etc.) allow us to
assume with out loss of generality that~$f = 1$. Thus, the problem
reduces to the following:

\subsubsection{}  One is given a~$p$-adic ring of integers~$A$, finite of
rank~$e$ over~$W$, and a free~$A$-algebra~$B$ of rank~$2$ with an
involution~$\ast$ whose fixed points are~$A$. Also given is a
perfect alternating pairing~$\langle \cdot, \cdot \rangle\colon
B\times B \rightarrow W$ such that for~$s\in B$ we have~$\langle
sr, r' \rangle = \langle r, s^\ast r' \rangle$. One needs to show
that every maximal isotropic~$B\otimes_W \fbar$ submodule
of~$B\otimes_W \fbar$ lifts to characteristic zero in the sense
previously described.

Note that~$B$ is either an integral domain that is a ramified
extension of~$A$ or isomorphic as an~$A$-algebra to~$A\oplus A$
with the involution being the permutation of coordinates. The
first case is immediate: We have~$B\otimes_W \fbar \cong
\fbar[t]/(t^{2e})$ and it has a unique submodule of rank~$e$
over~$\fbar$, viz. $(t^e)$. Since the Grassmann scheme~${\bf G'}$
always has characteristic zero geometric points and is projective,
a lift is provided by (any) characteristic zero point of~${\bf
G'}$.

In the second case we have~$B\otimes_W \fbar \cong \fbar[t]/(t^e)
\oplus \fbar[t]/(t^e)$. Every submodule of~$B\otimes_W \fbar$ of
rank~$e$ over~$\fbar$ is a direct sum~$(t^i)\oplus (t^{e-i})$.
Such submodules are automatically isotropic. We claim that the
submodule~$(t^i)$ of~$A[t]/(t^e)$ can be lifted to characteristic
zero, that such a lifting corresponds to a choice of~$e-i$
embeddings~$A \rightarrow \overline{\QQ}_p$ over $W$ and that each
lifting is isotropic when considered as a submodule of $B\otimes_W
W' = A\otimes_W W' \oplus A\otimes_W W'$, where~$W'$ is a ``big
enough" extension of~$W$. Indeed, every geometric point of the
appropriate Grassmann scheme, being proper over~$\Spec(W)$,
extends to an integral point (defined over a finite integral
extension~$W'/W$). Such a geometric point corresponds to a choice
of~$\left(\begin{smallmatrix} e \\ i
\end{smallmatrix} \right)$ embeddings~$A\rightarrow \overline{\QQ}_p$ over~$W$ and is isotropic
(cf.~\S\ref{subsubsection:3765} -- when we view~$A\otimes W$ as
a~$B$-submodule of~$B\otimes W$ via the first (or second)
component, it is isotropic). Moreover, since the submodule~$(t^i)$
is uniquely determined by its rank, every such integral point
indeed provides a lift of~$(t^i)$. It now easily follows
that~$(t^i)\oplus (t^{e-i})$ can be lifted
in~$\left(\begin{smallmatrix} e \\ i
\end{smallmatrix} \right)$ ways.

\subsection{Proof of Theorem~\ref{thm: embedding implies bad
reduction}} By an \emph{abelian scheme with CM} we mean in this
section a triple $(A/S, \lambda, \iota)$, consisting of a
principally polarized abelian scheme over~$S$ with an embedding of
rings $\iota\colon \calO_K \rightarrow \End_S(A)$ such that the
Rosati involution defined by $\lambda$ induces complex conjugation
on~$\calO_K$. We denote complex conjugation on ~$K$ by~$\ast$ and
let~$K^+$ be the totally real subfield of~$K$. As before, $W =
W(\fbar)$. The following lifting lemma, that holds for any CM
field~$K$ and whose proof is given in \S\S \ref{subsubsction local
triviality} - \ref{subsubsection end of proof}, is the key point.
\begin{lem}\label{lifting lemma} Let $(A, \lambda, \iota)$ be an abelian variety with CM over~$\fbar$
then $(A, \lambda, \iota)$ can be lifted to characteristic zero.
\end{lem}
\subsubsection{}\label{subsubsction local triviality}
Let~$S$ be a local artinian ring with residue field~$\fbar$.
Let~$(A', \lambda', \iota')$ be an abelian scheme over~$S$ with
CM. We claim that~$\HH^1_{\rm dR}(A'/S)$ is a free~$\calO_K
\otimes_\ZZ S$-module of rank~$1$. Since~$\HH^1_{\rm dR}(A'/S)$ is
a free~$S$-module of rank~$2g$, to verify that it is a
free~$\calO_K\otimes_\ZZ S$-module it is enough to prove that
modulo the maximal ideal of~$S$ (cf. \cite[Rmq. 2.8]{DP}), namely,
that $\HH^1_{\rm dR}(A'\otimes_S\fbar/\fbar)$ is a
free~$\calO_K\otimes_\ZZ \fbar$-module. This is \cite[Lem.
1.3]{Rapoport}. In fact, loc. cit. gives that $H^1_{\rm
crys}(A'\otimes_S\fbar/W)$ is a free~$\calO_K\otimes_\ZZ W$-module

\subsubsection{} The polarization~$\lambda$ induces a perfect
alternating pairing~$\langle \cdot, \cdot \rangle$ on the free
$\calO_K \otimes_\ZZ W$-module $H^1_{\rm crys}(A/W)$, which we
identify with $M_0:= \calO_K \otimes_\ZZ W$. This pairing induces
complex conjugation on~$\calO_K$ and reduces modulo~$p$ to the
pairing induced by~$\lambda$ on~$\HH^1_{\rm
dR}(A/\overline{\FF}_p)$. Moreover, there exists a finite
extension~$\Lambda$ of~$W$ such the Hodge filtration~$0
\rightarrow H^0(A, \Omega^1_{A/\overline{\FF}_p})\rightarrow
\HH^1_{\rm dR}(A/\overline{\FF}_p)$ can be lifted to~$M_0\otimes_W
\Lambda$. This follows from the discussion in~\S\ref{sec:
Grassmannians}. In fact, the results of that section show that
such a lift is uniquely determined by its generic point, a
subspace of $K \otimes_\QQ \overline{\QQ}_p =
\oplus_{\{\varphi\colon K \rightarrow\overline{\QQ}_p\}}
\overline{\QQ}_p(\varphi)$, consisting of a choice of one subspace
out of each pair $\overline{\QQ}_p(\varphi)\oplus
\overline{\QQ}_p(\varphi\circ \ast)$.

Recall that a CM type~$\Phi$ of~$K$ is a subset of~$\Hom(K, \CC)$
(or of $\Hom(K, \overline{\QQ}_p)$) that is disjoint from its
complex conjugate, equivalently, a subset that induces~$\Hom(K^+,
\CC)$ (or~$\Hom(K^+, \overline{\QQ}_p)$). A choice of lift of the
Hodge filtration provides us with CM type~$\Phi$. Let~$K^\ast$ be
the reflex field defined by~$\Phi$. We see that, in fact, a lift
of the Hodge filtration is defined over $\Lambda$, where $\Lambda$
is the compositum of $W$ with the valuation ring of the $p$-adic
reflex field associated to~$\Phi$.

\subsubsection{} Let~$V = \calO_K\otimes_\ZZ \CC$ - a complex vector
space on which~$\calO_K$ acts. Choose a~$\ZZ$-basis~$e_1, \dots,
e_{2g}$ for~$\calO_K$ and consider~$f_\Phi(\ux) = \det(\sum
e_ix_i, \Lie(A))$. This a polynomial in~$x_1, \dots, x_{2g}$ with
coefficients in~$\calO_{K^\ast}$ that depends only on~$\Phi$ and
determines it.

Let~${\bf M}\colon {\rm \underline{Sch}}_{\calO_{K^\ast}}
\rightarrow {\rm \underline{Sets}}\;$ be the functor from the
category of schemes over~$\calO_{K^\ast}$ to the category of sets
such that~${\bf M}(S)$ is the isomorphism classes of
triples~$(A/S, \lambda, \iota\colon \calO_K \injects \End_S(A))$,
where~$A/S$ is an abelian scheme with CM and~$\det(\sum e_ix_i,
\underline{\Lie}(A)) = f_\Phi(\ux)$. That is, the triple~$(A/S,
\lambda, \iota\colon \calO_K \injects \End_S(A))$ satisfies the
Kottwitz condition \cite[\S5]{Kottwitz Points} uniquely determined
by~$\Phi$.

For the given point~$x = (A, \lambda, \iota)\in {\bf M}(\fbar)$ we
consider the local deformation problem induced by~${\bf M}$. This
is the functor~${\bf N}$ from the category~${\bf C}_\Lambda$ of
local artinian~$\Lambda$-algebras with residue
field~$\overline{\FF}_p$ to the category \:\underline{Sets}\:
associating to a ring~$R$ in~${\bf C}_\Lambda$ those elements
of~${\bf M}(R)$ specializing to~$x$. We remark that the Kottwitz
condition is closed under specialization. It is thus fairly
standard that~${\bf N}$ is pro-represented by a complete
noetherian~$\Lambda$-algebra~$R^{\rm u}$; cf. \cite[\S2]{OortOslo}
and \cite[\S4]{Tim}.

\subsubsection{}\label{subsubsection end of proof} Let~${\bf G}\rightarrow \Spec(\Lambda)$ be the Grassmann variety
parameterizing for a scheme~$S\rightarrow \Spec(\Lambda)$ the set
of~$\calO_K$-invariant, isotropic, locally free, locally direct
summands~$\calO_S$-submodules of rank~$g$ of~$M_0\otimes_\Lambda
\calO_S$ (with the pairing coming from~$x$ as above) and
satisfying Kottwitz condition~$f_\Phi$ for a CM type~$\Phi$. (In
fact, one can deduce that ${\bf G} \cong \Spec(\Lambda)$ but we
don't need it here.) Let~$x$ be the point of~${\bf G}$
corresponding to~$H^0(A, \Omega^1_{A/\overline{\FF}_p})\rightarrow
\HH^1_{\rm dR}(A/\overline{\FF}_p)$.

Given the results of \S\ref{subsubsction local triviality}, the
theory of local models furnishes an isomorphism
$\calO^\wedge_{{\bf G}, x} \cong R^u$; cf. \cite[\S3]{DP},
\cite[Thm. 4.4.1]{Tim} -- the arguments easily extend to allow a
Kottwitz condition. We conclude therefore that there is a
triple~$(A, \lambda, \iota)$ lifting~$x$ defined over the~$p$-adic
field~$K_1 = \Lambda\otimes \QQ$. This concludes the proof of the
lemma.

\subsubsection{} Let~$K$ be a primitive quartic CM field. A solution of the
embedding problem for $p$ provides us with a triple $(A/\fbar,
\lambda, \iota)= (E_1\times E_2/\fbar, \lambda = \lambda_1\times
\lambda_2, \iota\colon \calO_K\rightarrow \End_{\fbar}(E_1\times
E_2))$. By Lemma~\ref{lifting lemma}, we may lift $(A/\fbar,
\lambda, \iota)$ to a triple $(A_0/\fbar, \lambda_0, \iota_0)$
defined over some~$p$-adic field~$K_1$ and so, by Lefschetz
principle, defined over~$\CC$. By the theory of complex
multiplication~$(A_0, \lambda_0, \iota_0)$ is defined over some
number field~$K_2$. Since the CM field~$K$ is primitive,~$A_0$ is
simple and principally polarized. By a theorem of Weil \cite{Weil}
the polarization is defined by a non singular projective genus~$2$
curve~$\scrC$ and it follows that~$A_0 \cong \Jac(\scrC)$ as
polarized abelian varieties. Furthermore,~$\scrC$ is defined over
a number field~$K_3$ (that is at most a quadratic extension
of~$K_2$). By passing to a finite extension~$L$ of~$K_3$, we get a
stable model.

\section{Applications}
\subsection{A general principle}
The following lemma is folklore and easy to prove:
\begin{lem}\label{lem: general principle} Let~$\pi\colon S\rightarrow R$ be a proper scheme over a Dedekind domain~$R$
with quotient field~$H$. Let~$\scrL \rightarrow S$ be a line
bundle on~$S$ and~$f, g\colon S\rightarrow \scrL$ sections.
Let~$x\in S(H')$ be a point, where~$H'$ is a finite field
extension of~$H$. Let~$u = (f/g)(x) \in H'$. Let~$\gerp$ be a
prime of~$R'$, the integral closure of~$R$ in~$H'$. Let~$\bar{x}$
be the~$R'$-point corresponding to~$x$. Then~$\val_\gerp(u) < 0$
implies that~$\bar{x}$ intersects the divisor of~$g$ in the fiber
of~$S$ over~$\gerp$.
\end{lem}
\begin{cor}\label{prime in denom implies bad reduction}
Let~$\scrA_{2}\rightarrow \Spec(\ZZ)$ be the moduli space of principally
polarized abelian surfaces and let
\begin{equation}\Theta(\tau) =
\frac{1}{2^{12}}\prod_{\begin{matrix}(\epsilon, \epsilon ') \\
{\rm even \; char.}\end{matrix}} \left(\Theta\left[
\begin{matrix} \epsilon \\ \epsilon ' \end{matrix}\right](0,
\tau)\right)^2.
\end{equation}
Let~$f$ be a Siegel modular form with~$q$-expansion~$\sum_{\nu}
a(\nu)q^{2\pi i \Tr(\:^t\nu  \tau)}$, where~$\nu$ runs over
$g\times g$ semi-integral, semi-definite symmetric matrices.
Assume that all the Fourier coefficients~$a(\nu)\in \calO_L$, the
ring of integers of a number field~$L$, and that the weight of~$f$
is of the form~$10 k$,~$k$ a positive integer.

Let~$\tau$ be a point on~${\rm Sp}_4(\ZZ)\backslash\gerH_2$
corresponding to a smooth genus~$2$ curve~$\scrC$ with CM by the
full ring of integers of a primitive CM field~$K$. Then
$(f/\Theta^k)(\tau)$ is an algebraic number lying in the
compositum~$H_{K^\ast} L$ of the Hilbert class field of~$K^\ast$
and~$L$. If a prime~$\gerp$ divides the denominator
of~$(f/\Theta^k)(\tau)$ then~$\scrC$ has bad reduction
modulo~$\gerp$.
\end{cor}
\begin{proof} The argument is essentially that of \cite[\S4.4]{De Shalit
Goren}: Igusa \cite{Igusa} proved that~$\Theta$ is a modular form
on~${\rm Sp}_4(\ZZ)\backslash\gerH_2$ (see \cite[Thm. 3]{Igusa},
$\Theta$ is denoted there~$\chi_{10}$). It is well known to have
weight~$10$ and a computation shows that its Fourier coefficients
are integers and have g.c.d.~$1$. The~$q$-expansion
principle~\cite[Ch. V, Prop. 1.8]{FaltingsChai} shows that~$f$ and
$\Theta^k$ are sections of a suitable line bundle of the moduli
scheme~$\scrA_2\otimes_\ZZ \calO_L$. The value
$(f/\Theta^k)(\tau)$ lies in~$H_{K^\ast}L$ by the theory of
complex multiplication.

It is classical that the divisor of~$\Theta$ over~$\CC$,
say~$D_{\rm gen}$, is the locus of the reducible polarized abelian
surfaces -- those that are a product of elliptic curves with the
product polarization. The Zariski closure~$D_{\rm gen}^{\rm cl}$
of~$D_{\rm gen}$ in~$\scrA_2$ is contained in the divisor~$D_{\rm
arith}$ of~$\Theta$, viewed as a section of a line bundle
over~$\scrA_2$, and therefore~$D_{\rm gen}^{\rm cl} = D_{\rm
arith}$, because by the~$q$-expansion principle~$D_{\rm arith}$
has no ``vertical components". Since~$D_{\rm gen}^{\rm cl}$ also
parameterizes reducible polarized abelian surfaces, it follows
that~$D_{\rm arith}$ parameterizes reducible polarized abelian
surfaces. (Furthermore, it is easy to see by lifting that every
reducible polarized abelian surface is parameterized by~~$D_{\rm
arith}$.) The Corollary thus follows from Lemma~\ref{lem: general
principle}.
\end{proof}
\begin{cor}\label{bounding corollary}$(f/\Theta^k)(\tau)$ is an~$S$-integer, where~$S$ is
the set of primes of lying over rational primes~$p$ less than
$16\cdot  d^2 (\Tr(r))^2$ and such that~$p$ decomposes in a
certain fashion in a normal closure of~$K$ as imposed by
superspecial reduction~\cite[Thms. 1, 2]{Goren Reduction} (for
example, if~$K$ is a cyclic Galois extension then~$p$ is either
ramified or decomposes as~$\gerp_1\gerp_2$ in~$K$).
\end{cor}

\subsection{Class invariants} Igusa \cite[p. 620]{Igusa AVM}
defined invariants $A(u), B(u), C(u), D(u)$ of a sextic $u_0X^6 +
u_1X^5 + \cdots + u_6$, with roots $\alpha_1, \dots, \alpha_6$, as
certain symmetric functions of the roots. For example, $D(u) =
u_0^{10} \prod_{i<j}(\alpha_i - \alpha_j)^2$ is the discriminant.
Igusa also proved that if $k$ is a field of characteristic
different from $2$, the complement of $D = 0$ in ${\rm
\text{Proj}\;}k[A, B, C, D]$, where $A, B, C, D$ are of weights
$2, 4, 6, 10$ respectively, is the coarse moduli space for
hyperelliptic curves of genus $2$. Moreover, the ring of rational
functions is generated by the ``absolute invariants" $B/A^2,
C/A^3, D/A^5$ (see \cite[p. 177]{Igusa1}, \cite[p. 638]{Igusa
AVM}). One can choose other generators of course, and for our
purposes it makes sense to choose generators with denominator a
power of $D$. Choose then as in \cite[p. 313]{van Wamelen} the
generators
\[i_1 = A^5/D ,\quad i_2 = A^3B/D, \quad i_3 = A^2C/D.\]
One should note though that these invariants are not known
a-priori to be valid in characteristic~$2$, since Weierstrass
points ``do not reduce well" modulo $2$. The invariants $i_n$ can
be expressed in terms of Siegel modular forms thus:
\[ i_1 = 2\cdot 3^5 \chi_{10}^{-6}\chi_{12}^5, \quad i_2 =
2^{-3}\cdot 3^{3} \psi_4 \chi_{10}^{-4}\chi_{12}^3, \quad i_3 =
2^{-5}\cdot 3 \psi_6 \chi_{10}^{-3}\chi_{12}^2 + 2^{2}\cdot 3
\psi_4 \chi_{10}^{-4}\chi_{12}^3.\] See \cite[pp. 189,
195]{Igusa1} for the definitions; $\psi_i$ are Eisenstein series
of weight $i$, $-2^{2}\chi_{10}$ is our $\Theta$.

Another interesting approach to the definition of invariants is
the following: Let $I_2 = h_{12}/h_{10}, $ $I_4 = h_4, \, I_6 =
h_{16}/h_{10}, \, I_{10} = h_{10}$ be the modular forms of weight
$2, 4, 6, 10$, respectively, as in \cite{Lauter Class Invariants}.
The appeal of this construction is that each $h_n$ is a simple
polynomial expression in Riemann theta functions with integral
even characteristics $\left[
\begin{smallmatrix} \epsilon \\ \epsilon ' \end{smallmatrix}\right]$; for example,
$h_4 = \sum_{10} (\Theta\left[
\begin{smallmatrix} \epsilon \\ \epsilon ' \end{smallmatrix}\right](0, \tau))^4 $, $h_{10} = 2^{12}\Theta$. It is not hard to
prove that the g.c.d. of the Fourier coefficients of $\Theta\left[
\begin{smallmatrix} \epsilon \\ \epsilon ' \end{smallmatrix}\right](0, \tau)$,
for $\left[\begin{smallmatrix} \epsilon \\ \epsilon '
\end{smallmatrix}\right]$ an integral even characteristic, is $1$
if $\epsilon\in \ZZ^2$ (that happens for $4$ even characteristics)
and $2$ if $\epsilon\not\in \ZZ^2$ (that happens for $6$ even
characteristics). Using that and writing $I_n = \ast/\Theta $, one
finds that the numerator of $I_n$ has an integral Fourier
expansion. One then lets
\[\gerj_1: = I_2^5/2^{-12} I_{10},
\quad \gerj_2: = I_2^3I_4/2^{-12}I_{10}, \quad \gerj_3: =
I_2^2I_6/2^{-12}I_{10}. \] These are modular functions of the form
$f/\Theta^k$, such that the numerator has integral Fourier
coefficients. Slightly modifying the definition of \cite{Lauter
Class Invariants} (there one uses $j_i:=2^{-12}\gerj_i$), we put
\begin{equation}\label{class polynomials 2} H_i(X) = \prod_\tau(X -
\gerj_i(\tau)), \quad i = 1, 2, 3,
\end{equation}
where the product is taken over all $\tau\in {\rm Sp}(4,
\ZZ)\backslash \gerH_2$ such that the associated abelian variety
has CM by $\calO_K$ (thus all polarizations and CM types appear).
We remark that $j_1 = i_1, j_2 = i_2$; this can be verified using
the formulas given in \cite[p. 848]{Igusa ModForms}.

The polynomials appearing in Equation (\ref{class polynomials 2})
have rational coefficients that are symmetric functions in modular
invariants, viz. the values of the functions $\gerj_i$ associated
to CM points. As such, it is natural to ask for the prime
factorization of these coefficients. For example, the results
of~\cite{Gross Zagier} give the factorization of the discriminant
of the Hilbert class polynomial in the case of imaginary quadratic
fields and so provide a bound on the primes which can appear. In
\cite{Lauter Class Invariants}, it was conjectured that primes
dividing the denominators of the coefficients of $H_i(X)$ are
bounded by the discriminant of $K$ (note that the only difference
between the current definition and loc. cit. is powers of $2$). We
deduce from the preceding results the following:
\begin{cor} The coefficients of the rational polynomials $H_i(X)$
are $S$-integers where $S$ is the set of primes smaller than
$16\cdot d^2 (\Tr(r))^2$ and satisfying a certain decomposition
property in a normal closure of $K$ as imposed by superspecial
reduction.
\end{cor}
We provide some numerical data in \S\S\ref{app:class invariants} -
\ref{app: curves with bad reduction}.
\begin{rmk} Theorem~\ref{thm: embedding implies bad reduction}
gives a partial converse to this corollary.
\end{rmk}

\subsection{Units} Let $K$ be a primitive quartic CM field as before.
In \cite{De Shalit Goren}, De Shalit and the first named author
constructed class invariants $u(\Phi; \gera)$, $u(\Phi; \gera,
\gerb)$ associated to certain ideals of $K$ and a CM type $\Phi$.
The construction essentially involves the evaluation of $\Theta$
at various CM points associated to $K$. Though the construction is
general, we recall it only for the $u(\Phi;\gera)$ and under very
special conditions. For the general case, refer to loc. cit.
\begin{exa} Assume that $K$ is a cyclic CM field with odd class
number $h_K$, $h_{K^+} = 1$. Let $\Phi$ be a CM type of $K$ and
assume that the different ideal $\calD_{K/\QQ} = (\delta)$ with
$\overline{\delta} = -\delta$ and $\Ima(\varphi(\delta)) > 0$ for
$\varphi\in \Phi$. Let $\gera$ be a fractional ideal of $\calO_K$
and choose $a\in K^+, a\gg 0$ such that $\gera\overline{\gera} =
(a)$. The form $\langle f, g \rangle =
\Tr_{K/\QQ}(\overline{f}g/a\delta)$ induces a principal
polarization on $\CC^2/\Phi(\gera)$. Write the lattice
$\Phi(\gera)$ as spanned by the symplectic basis formed by the
columns of $(\omega_1 \omega_2)$ and consider
$\Delta(\Phi(\gera)):=
\det(\omega_2)^{-10}\Theta(\omega_2^{-1}\omega_1)$. It depends
only on $\Phi, \gera$ and not on $a$. One then lets
\begin{equation} u(\Phi;\gera) =
\frac{\Delta(\Phi(\gera^{-1}))}{\Delta(\Phi(\calO_K))}, \qquad
u(\Phi; \gera, \gerb) = \frac{u(\Phi; \gera\gerb)}{u(\Phi;\gera)
u(\Phi;\gerb)}.
\end{equation}
See \cite[\S1.3]{De Shalit Goren} for remarkable properties of
these invariants. In particular, if $h_K$ is a prime different
from $5$ the group generated by the $u(\Phi; \gera, \gerb)$ in
$H_K^\times$ has rank $h_K - 1$. The following corollary holds in
general.
\end{exa}
\begin{cor} The invariants $u(\Phi; \gera, \gerb)$ are $S$-units
for $S$ the set of primes of $H_{K^\ast}$ that lie over rational
primes $p$ smaller than $16\cdot d^2 (\Tr(r))^2$ such that $p$
decomposes in a certain fashion in a normal closure of $K$ as
imposed by superspecial reduction.
\end{cor}

\section{Appendix: Numerical data}

\subsection{Class invariants}\label{app:class invariants}

Let $K=\QQ[x]/(x^4+50x^2+93)$ be the non-normal quartic CM field
of class number 4 generated by $i\sqrt{25+2\sqrt{133}}$ over its
totally real subfield $K_0 = \QQ(\sqrt{133})$.  The field
discriminant of $K$ is $d = 3 \cdot 31 \cdot 133^2$.  The reflex
field of $K$ is the quartic CM field
$K^*=\QQ[x]/(x^4+100x^2+2128)$, and it also has class number 4.
The first class polynomial $H_1(X)$ for $K$ is:

\begin{flushleft} \begin{tiny}
$H_1(X) = X^8 + (1041267141265383834470066376076324878559188075169214718620644562742786358262348341096508075804599$ \\
$14396998445377520399898369594754800397259854505319134527276855203669924179304812893406433677385727145618982259878$ \\
$235238106861354975708604204456451607420809338475387838838399714721792 X^7 - 12149397793963178112627821022620892$\\
$45404136087671598296953472940404425708214322876976589037187761355014806902217382933228421664445385161486770241395$ \\
$86216425653442970815869576705927409199788135906841889599575599439952333875289149403700489027842584764660920021436$ \\
$3829932461236811798235188297728 X^6 + 73839049971416349806076414527300316633378897475063812302053799419506624$\\
$92518129165964459009432807362787172373507594090454704865662995335680600282096554468522204003852710171918281867104$ \\
$22881084398565227843015748927435317832813123246139494220230411879998467816549366023387634647423494310078523116827$ \\
$4432 X^5 + 13572072584610002627351064096562440800687789906069625620954572275815407089929054924847221405346617$\\
$81303348821546686327912447263712909341876160159292258061730849139344294186108556484931948397024260104673954036849$\\
$458602462351542470926873770897062713438092865784651432248872830418650819961887941339751759032979816448 X^4 - $ \\
$20645419071899094062319060950215185132822523834337866324899816490099930519115325291745691785099541742412841978118$\\
$50659093101670013797324826062128177308487589839435711979382829163139412065631102085098562891221391707621189640462$\\
$6733170257389221578842790790536175427848548886553551634053100728321789568362872897268981614772224 X^3 - 129312493$\\
$7479769249061778818022624358655510743462873298439586515036372146763574796362009513855612419032661986690064195$\\
$42745542696873991073642021736344937793632645966542526675266553096755940524528853447553985574064564786043667420743$\\
$99008810426655744434090152803879394309712176142971045663487522627933115839853095635881434146324611072 X^2 - $\\
$175163105404286285404558242388178232401134785649487110037435573206470551113657779259601408409259528822541237329511$\\
$265595839919081218671703702189931407432529268004641208854959298642809557864657052618947314472582198458807634678798$\\
$21263464159806719108635491851367067315658512258914788283117920745406313054218903105943266111760957022193902944256 X$\\
$+ 1604398176398171965083954530143587467864124228892266410893720503328183568724268331605647208243792414122$\\
$982980943620961424734744552992794289167983185212334500945230977070579478907220448256281603461621068539088549163272$\\
$6437672763292720798162365206307314646594632579281296420073078445474401380548434319635226631450753809199086032841732$\\
$365137199566422016)/16549715179319233558563819433380761888575722585072935139257013386884308496619078534362406415399$\\
$3060472421943063977491964645435818434816422816448520354952124186248583885336013832737212668433932405635280914057910$\\
$76698802759845571321099052922204371789683841$\\
\end{tiny}
\end{flushleft}

The other two class polynomials are not given here, as they have
the same primes in their denominators.  The polynomials were
computed using PARI with 1000 digits of precision in about 8 hours
each on an Intel Pentium 4, 2.2GHZ, 512MB memory. The denominator
factors as: $7^{48} \cdot 11^{72} \cdot 19^{24} \cdot 23^{12}
\cdot 29^{12} \cdot 83^{12} \cdot 89^{12} \cdot 167^{12}.$ Note
that for the first invariant $\Theta$ appears to the sixth power
in the denominator, which agrees with the fact that all powers are
a multiple of $6$.

\subsection{Curves with bad reduction}\label{app: curves with bad reduction}
To illustrate the theory we give an example of a CM field $K$ and
two genus 2 curves over $\QQ$ with CM by $K$.  We list their
invariants, and verify that they have bad reduction at the primes
in the denominators of the invariants. In~\cite{van Wamelen}, van
Wamelen gives a complete list of all isomorphism classes of genus
2 CM curves defined over the rationals along with their Igusa
invariants.  For example, for the cyclic CM field
$K=\QQ(i\sqrt{13-3\sqrt{13}})$ of class number 2, there are two
non-isomorphic genus 2 curves defined over $\QQ$.

The curve with invariants equal to $i_1=\frac{2\cdot 11^5 \cdot
53^5 \cdot 6719^5 \cdot 30113^5 }{3^7 \cdot 23^{12} \cdot
131^{12}}$, $i_2 =\frac{2\cdot 5 \cdot 11^3 \cdot 53^3 \cdot
6719^3 \cdot 7229 \cdot 30113^3 }{3^3 \cdot 23^{8} \cdot
131^{8}}$, $i_3 =\frac{2\cdot  11^2 \cdot 19 \cdot 53^2 \cdot
6719^2 \cdot 30113^2 \cdot 237589628623651}{3^4 \cdot 23^{8} \cdot
131^{8}},$ has an affine model
$$y^2 = -70399443x^6+36128207x^5+262678342x^4-48855486x^3-112312588x^2+36312676x. $$
The reduction of a genus 2 curve at a prime can be calculated
using~\cite[Thm 1, p. 204]{Liu}. For these examples we actually
calculated the reduction using the genus 2 reduction program
written by Liu. The output of the program shows that at the primes
$p=2,\,3\,,\,23,\,131,$ the curve has potential stable reduction
equal to the union of two supersingular elliptic curves $E_1$ and
$E_2$ intersecting transversally at one point.

The second curve has invariants equal to $i_1=\frac{2\cdot 7^{10}
\cdot 11^5 \cdot 21059^5}{3^7 \cdot 23^{12}}$, $i_2 =\frac{2\cdot
5 \cdot 7^7 \cdot 11^3 \cdot 8387 \cdot 21059^3}{3^3 \cdot
23^{8}},$ $i_3 =\frac{2\cdot 7^6 \cdot 11^2 \cdot 21059^2 \cdot
71347 \cdot 739363}{3^4 \cdot 23^{8}},$ and has an affine model
$$y^2 = -243x^6+2223x^5-1566x^4-19012x^3+903x^2+19041x-5882. $$
In this case, the output of the genus 2 reduction program shows
that at the primes $p=2,\,3\,,\,23,$ the curve has potential
stable reduction equal to the union of two supersingular elliptic
curves $E_1$ and $E_2$ intersecting transversally at one point.

The reader may have noticed that $2$ does not appear in the
denominator of the invariants. This is not due to the invariants
$i_n$ being divisible by $2$. It is an artifact of cancellation
between ``values of the numerator and the denominator" and
explains in which sense Theorem~\ref{thm: embedding implies bad
reduction} may fail to provide a converse to Corollary~\ref{prime
in denom implies bad reduction}. In fact, bad reduction of CM
curves modulo primes over 2 turns out to be prevalent. According
to~\cite{IKO}, there is no smooth superspecial curve in
characteristic 2. On the other hand, using complex multiplication,
one can prove (e.g. for cyclic CM fields K and primes decomposing
as $(p) = \gerp_1 \gerp_2$ or $(p) = \gerp_1^2$) superspecial
reduction of principally polarized abelian surfaces with CM by K
(cf.~\cite{Goren Reduction}). This implies for $p=2$ bad reduction
of the corresponding curve.

\subsection{Isogenies between two supersingular curves}\label{data:isogenies}Let~$p$ be a prime,~$h$ the
class number of~$B_{p, \infty}$ and~$N(p)$ the minimal integer for
which there exists an isogeny of degree less or equal to~$N(p)$
between any two supersingular elliptic curves over
$\overline{\FF}_p$. Because of running time and memory
restrictions we did only sample calculations. For~$p$ = 10007, the
Total Computation Time was 22688.710 seconds, Total Memory Usage
was 1213.97MB. The program ran on an Intel Pentium 4, 2.53 GHz, 1
GB memory.

\vspace{0.5cm}

{\small
\begin{tabular}{|l|l|l|l|l|}\hline
$p$ &      ~$h$ &      ~$[\sqrt{p}]$  &~$N$ &      ~$N/\sqrt{p}$ \\
\hline\hline 101   &  9 &      10 & 6 &      0.600
\\\hline
211 & 18 & 15& 9& 0.600
\\\hline
307 & 26   &   18  & 11 &     0.611
\\\hline
401 & 34 &     20 &12& 0.600
\\\hline
503 & 43 &22& 15  &    0.682
\\\hline
601 & 50 &25 &14 &0.560
\\\hline
701 & 59& 26& 17& 0.654
\\\hline
809 &68&  28 &18& 0.643
\\\hline
907 & 76& 30 &19 &0.633
\\\hline
1009& 84 &32& 20& 0.625
\\\hline
2003& 168& 45& 30 &0.667
\\\hline
3001 & 250 &55 &34 &0.618
\\\hline
4001 &334 &63& 44& 0.698
\\\hline
5003 & 418 &71& 46& 0.648
\\\hline
6007 &501& 78& 51& 0.654
\\\hline
7001 & 584 &84& 56& 0.667
\\\hline
8009 &668& 89 &60 &0.674
\\\hline
9001 &750 &95 &59 &0.621
\\\hline
10007 &835 &100 &70 &0.700
\\\hline \end{tabular}
}

\

\

\id {\emph{Acknowledgments: We would like to thank Payman Kassaei
for valuable discussions and Ehud DeShalit for some useful
comments. The first named author's research was partially
supported by NSERC; he would like to thank Microsoft Research for
its hospitality during a visit where this project took shape.  The
second named author thanks Tonghai Yang for many stimulating
discussions and McGill University for its hospitality.}}

\end{document}